\documentclass[letterpaper, 10 pt, conference]{ieeeconf}  

\IEEEoverridecommandlockouts                              

\usepackage{color}
\usepackage{graphics} 
\usepackage{epsfig} 
\usepackage{subfigure}
\usepackage{ulem}
\usepackage{times} 
\usepackage{amsmath} 
\usepackage{amssymb}  
\usepackage{tikz}

\newtheorem{theorem}{Theorem}[section]
\newtheorem{proposition}{Proposition}[section]
\newtheorem{definition}{Definition}[section]

\newtheorem{remarkth}[definition]{Remark}
\newenvironment{remark}{\begin{remarkth}\upshape}{\end{remarkth}}

\usepackage{amssymb}

\newcommand{\proa}{A^*G \mbox{$\;$}_{\tau^*} \kern-3pt\times_\alpha
G \mbox{$\;$}_\beta \kern-3pt\times_{\tau^*} A^*G}

\newcommand{\ad}{\mbox{ad}}

\title{\LARGE \bf
Optimal Control of Left-Invariant Multi-Agent Systems with Asymmetric Formation Constraints}

\author{Leonardo J. Colombo and Dimos V. Dimarogonas
\thanks{L.J. Colombo (colombo2@kth.se) and D.V. Dimarogonas (dimos@kth.se) are with ACCESS Linnaeus Center, Department of Automatic Control, KTH, Royal Institute of Technology, SE-100 44, Stockholm, Sweden. This work was supported by the Swedish Research Council (VR), Knut och Alice Wallenberg foundation (KAW), the H2020 Project Co4Robots, the H2020 ERC Starting Grant BUCOPHSYS and MINECO (Spain) grant MTM2016-76072-P.}}%
\begin{document}

\maketitle
\thispagestyle{empty}
\pagestyle{empty}

\begin{abstract}
In this work we study an optimal control problem for a multi-agent system modeled by an undirected formation graph with nodes describing the kinematics of each agent, given by a left invariant control system on a Lie group. The agents should avoid collision between them in the workspace. Such a task is done by introducing some potential functions into the cost functional for the optimal control problem, corresponding to fictitious forces, induced by the formation constraint among agents, that break the symmetry of the individual agents and the cost functions, and rendering the optimal control problem partially invariant by a Lie group of symmetries. Reduced necessary conditions for the existence of normal extremals are obtained using techniques of variational calculus on manifolds. As an application we study an optimal control problem for multiple unicycles. \end{abstract}

\section{Introduction}

Multi-agent systems \cite{MMbook} have gained a lot of attention in the last decades due to the advances in communication, robotics and cooperative control of spacecrafts, robotic manipulators as well as unmanned aerial and underwater vehicles. Symmetries in optimal control have been studied by several authors in the last decades \cite{Bl}, \cite{Krishnaprasad}, \cite{Leonard1} and the reduction by symmetries in optimal control problems have been a very active area of research for their applications in robotics, aerospace engineering and locomotion among others (see e.g., \cite{CoMdD}, \cite{KoonMarsden} and references therein).

Different approaches to formation control of multi-agent systems can be identified, e.g. as leader-follower \cite{Leonard2}, behavior-based \cite{balch} and rigid body type formations \cite{pedro1}, \cite{pedro2}, \cite{tabuada}, \cite{chris}. We build on the last category by studying optimal control of formation problems for systems whose configurations evolves on a Lie group and including in our analysis rigid body dynamics.

This work follows the research lines started in \cite{Krishnaprasad}  for optimal control of left-invariant systems and \cite{JK} for coordination control of multiple left invariant agents, respectively, and also builds in previous developments for reduction of optimal control \cite{BlCoGuOh}, \cite{KoonMarsden}, \cite{Krishnaprasad} by studying optimal control problems for multi-agent formations whose dynamics evolves on a Lie group of symmetries and the kinematic of each agent is given by a left-invariant system. 

The problem studied in this work consists on finding the absolute configurations and control inputs for each agent, obeying the corresponding kinematics equations given by a left-invariant control system, as well as satisfying the formation constraints and minimizing the energy of the agents in the formation.

One of the aims of this contribution, further than only solving the proposed optimal control problem, consists in introducing a new theoretical approach for the optimal coordinated motion of multi-agent systems with heterogeneous agents using variational principles, as is usual for a single agent whose dynamics is described by a mechanical system. The differential equations obtained from the variational principle and representing necessary conditions for optimality gives rise to a set of equations that can not be solved directly with a numerical solver. We also propose a splitting in the equations to render such a system solvable.

In this work each agent is described by a drift-free kinematic control system on a Lie group and agents should satisfy the formation constraints to avoid collision with each other in the workspace. For this task, we introduce some potential functions corresponding to fictitious forces, induced by the formation constraints, into the cost functional for the optimal control problem. Such potential functions are not invariants under the group of symmetries of the agents and therefore they break the symmetry of individual agents in the optimal control problem. The reduction of necessary conditions for the existence of normal extremals in the problem is described via Euler-Poincar\'e equations arising from the variational analysis.
 
As an application we study an energy-minimum problem for three unicycles and characterize the exact solution for the control law of one of the agents. To the best of authors' knowledge, this is one of the first attempts where a formation constraint for a coordination motion of unicycles is expressed in absolute configurations on the Lie group $SE(2)$ (a different approach for relative configurations has been studied in \cite{JK}), allowing to explore more the Lie group framework in formation problems with non-compact configuration spaces. This approach can be seen as a complement to the related literature for formation problems on Lie groups is the case of agents evolving on the Lie group of rotations $SO(3)$ where the constraint is written as the geodesic distance between two points (since $SO(3)$ is compact and therefore a complete manifold such a distance is well defined) and the use of Rodrigues' formula allows the use of trackable formation constraints.

The structure of the paper is as follows: Section II introduces reduction by symmetries and Euler-Poincar\'e equations. Section III introduces the  left-invariant kinematic multi-agent control system and the formulation of the optimal control problem for multiple agents. In Section IV we study Euler-Poincar\'e  reduction of necessary conditions by using a variational framework and splitting the dynamics to find a solvable system of equations. An application to the optimal control problem of three unicycles is studied in Section V. We conclude the work by commenting directions of future research in Section VI. 
\section{Preliminaries: Mechanics on manifolds}\label{section2}

In this section we introduce the review material we will use along the work. For a further covering of the topics see \cite{Bl} (Chapter $3$), \cite{Holmbook} (Chapters $6$-$7$), and \cite{Sastrybook} (Appendix A).

\subsection{Mechanics on manifolds}
Let $Q$ be the configuration space of a mechanical system, a differentiable manifold of dimension $n$ with local coordinates $q=(q^1,\ldots,q^n)$. Let $TQ$ be the cotangent bundle of $Q$, locally described by positions and velocities for the system $v_q=(q^1,\ldots, q^n,\dot{q}^{1},\ldots,\dot{q}^{n})\in TQ$ with $\hbox{dim}(TQ)=2n$. 

The dynamics of the mechanical system is determined by a Lagrangian function $L:TQ\to\mathbb{R}$ given by $L(q,\dot{q})=K(q,\dot{q})-V(q)$ where $K:TQ\to\mathbb{R}$ is the kinetic energy and $V:Q\to\mathbb{R}$ the potential energy. The corresponding equations of motion describing the dynamics of the system are given by the Euler-Lagrange equations $\displaystyle{\frac{d}{dt}\left(\frac{\partial L}{\partial\dot{q}^{i}}\right)=\frac{\partial L}{\partial q^{i}}}$, $i=1,\ldots,n;$ a system of $n$ second-order differential equations.
\subsection{Mechanics on Lie groups}
If the configuration space is a Lie group and the system has a symmetry, one can take advantage of it to reduce the degrees of freedom of the system and work on a lower dimensional system reducing computational cost and avoiding singularities by working on a coordinate free framework in the associated Lie algebra of a Lie group.
\begin{definition} A \textit{Lie group} is a smooth manifold $G$ that is a group and for which the operations of multiplication $(g,h)\mapsto gh$ for $g,h\in G$ and inversion, $g\mapsto g^{-1}$, are smooth.

\end{definition}
\begin{definition}
A \textit{symmetry} of a function $F:G\to\mathbb{R}$ is a map $\phi:G\to G$ such that $F\circ\phi=F$.
\end{definition}
\begin{definition} Let $G$ be a Lie group with identity element $\bar{e}\in G$. A \textit{left-action} of $G$ on a manifold $Q$ is a smooth mapping $\Phi:G\times Q\to Q$ such that $\Phi(\overline{e},q)=q$ $\forall q\in Q$, $\Phi(g,\Phi(h,q))=\Phi(gh,q)$ $\forall g,h\in G, q\in Q$ and for every $g\in G$, $\Phi_g:Q\to Q$ defined by $\Phi_{g}(q):=\Phi(g,q)$ is a diffeomorphism. \end{definition}

We often use the notation $gq:=\Phi_{g}(q)=\Phi(g,q)$ and say that $g$ acts on $q$. All actions of Lie groups will be assumed to be smooth.
Let $G$ be a finite dimensional Lie group and $\mathfrak{g}$ will denote the Lie algebra associated to $G$ defined as $\mathfrak{g}:=T_{\overline{e}}G$, the tangent space at the identity $\overline{e}\in G$.  Let $L_{g}:G\to G$ be the left translation of the element $g\in G$ given by $L_{g}(h)=gh$ for $h\in G$. $L_g$ is a diffeomorphism on $G$ and a left-action from $G$ to $G$ \cite{Holmbook}. Their tangent map  (i.e, the linearization or tangent lift) is denoted by $T_{h}L_{g}:T_{h}G\to T_{gh}G$. Similarly, the cotangent map (cotangent lift) is denoted by $T_{h}^{*}L_{g}:T^{*}_{h}G\to T^{*}_{gh}G$. It is well known that the tangent and cotangent lift are actions (see \cite{Holmbook}, Chapter $6$).

Let $X:G\to TG$ be a vector field on $G$. The set $\mathfrak{X}(G)$ denotes the set of all vector fields on $G$. The tangent map $T_{\overline{e}}L_g$ shifts vectors based at $\overline{e}$ to vectors based at $g\in G$. By doing this operation for every $g\in G$ we define a vector field as $X_{\xi}^{g}:=T_{\overline{e}}L_g(\xi)$ for $\xi:=X(\overline{e})\in T_{\overline{e}}G$. A vector field $X\in\mathfrak{X}(G)$ is called \textit{left-invariant} if $T_{h}L_g(X(h))=X(L_g(h))=X(gh)$ for all $g,h\in G$. In particular for $h=\overline{e}$ this means that a vector field $X$ is left-invariant if $\dot{g}=X(g)=T_{\overline{e}}L_{g}\xi$ for $\xi=X(\overline{e})\in\mathfrak{g}$. Note that if $X$ is a left invariant vector field, then $\xi=X(\overline{e})=T_{g}L_{g^{-1}}\dot{g}$.   

Let $\Phi_g:Q\to Q$ for any $g\in G$ a left action on $G$; a function $f:Q\to\mathbb{R}$ is said to be \textit{invariant} under the action $\Phi_g$, if $f\circ\Phi_g=f$, for any $g\in G$ (that is, $\Phi_g$ is a symmetry of $f$). The \textit{Adjoint action}, denoted $\hbox{Ad}_{g}:\mathfrak{g}\to\mathfrak{g}$ is defined by $\hbox{Ad}_{g}\chi:=g\chi g^{-1}$ where $\chi\in\mathfrak{g}$. Note that this action represents a change of basis on the Lie algebra.


If we assume that the Lagrangian $L\colon TG\to\mathbb{R}$ is $G$-invariant under the tangent lift of left translations, that is $L\circ T_{g}L_{g^{-1}}=L$ with $g\in G$, then it is possible to obtain a reduced Lagrangian $\ell\colon\mathfrak{g}\to\mathbb{R}$, where $$\ell(\xi) = L(g^{-1}g,T_{g}L_{g^{-1}}(\dot{g}))= L(\overline{e},\xi).$$The reduced Euler--Lagrange equations, that is, the Euler--Poincar{\'e} equations (see, e.g., \cite{Bl}, \cite{Holmbook}), are given by the system of $n$ first order ode's \begin{align}
\frac{d}{dt}\frac{\partial\ell}{\partial\xi} = \ad^{*}_{\xi}\frac{\partial\ell}{\partial\xi}.\label{eq_ep_intro}
\end{align} where $\ad^{*}:\mathfrak{g}\times\mathfrak{g}^{*}\to\mathfrak{g}^{*}$, $(\xi,\mu)\mapsto\ad^{*}_{\xi}\mu$ is the \textit{co-adjoint operator} defined by $\langle\ad_{\xi}^{*}\mu,\eta\rangle=\langle\mu,\ad_{\xi}\eta\rangle$ for all $\eta\in\mathfrak{g}$ with $\ad:\mathfrak{g}\times\mathfrak{g}\to\mathfrak{g}$ the \textit{adjoint operator} given by $\ad_{\xi}\eta:=[\xi,\eta]$, where $[\cdot,\cdot]$ denotes the Lie bracket of vector fields on the Lie algebra $\mathfrak{g}$, and where $\langle\cdot,\cdot\rangle:\mathfrak{g}^{*}\times\mathfrak{g}\to\mathbb{R}$ denotes the so-called \textit{natural pairing} between vectors and co-vectors defined by $\langle\alpha,\beta\rangle:=\alpha\cdot\beta$ for $\alpha\in\mathfrak{g}^{*}$ and $\beta\in\mathfrak{g}$ where $\alpha$ is understood as a row vector and $\beta$ a column vector. For matrix Lie algebras $\langle\alpha,\beta\rangle=\alpha^{T}\beta$ (see \cite{Holmbook}, Section $2.3$ pp.$72$ for details).

Using this pairing between vectors and co-vectors, one can write a useful relation between the tangent and cotangent lifts \begin{equation}\label{relation-cltl}\langle\alpha,T_{h}L_{g}(\beta)\rangle=\langle T^{*}_{h}L_{g}(\alpha),\beta\rangle\end{equation} for $g,h\in G$, $\alpha\in\mathfrak{g}^{*}$ and $\beta\in\mathfrak{g}$.

The Euler--Poincar{\'e} equations together with the reconstruction equation $\xi = T_{g}L_{g^{-1}}(\dot{g})$ are equivalent to the Euler--Lagrange equations on $G$.

\section{Left-invariant kinematic multi-agent control system and problem formulation}
\subsection{Left-invariant kinematic multi-agent control system.}
Let $G_i$ be $r$ Lie groups of dimension  $n$ describing the configuration of $r$ heterogenous agents and $\mathfrak{g}_i:=T_{\overline{e}_{i}}G_{i}$ their corresponding Lie algebras $i=1,\ldots,r$ with $g_i(t)\in G_i$ describing the evolution of agent $i$ at time $t$.

In the problem studied in this work the configuration space of each agent has the same Lie group structure. Note that the same configuration does not mean the same agent. For instance, each agent can have different mass and inertia values, and therefore agents are heterogeneous agents. 

Along this work, we assume that a multi-agent control system is modeled by an undirected (bidirectional) formation graph $\mathcal{G}=(\mathcal{V},\mathcal{E}, \mathcal{C})$, describing the kinematics of each agent given by $r$ left invariant control systems on $G_i$ with $i=1,\ldots,r$ together with the formation constraints. 

Here $\mathcal{V}$ denotes the set of vertices of the graph representing the communication topology in the multi agent system where each vertex $v_i\in\mathcal{V}$ is a left invariant control system, that is, the kinematics of each agent is determined by  \begin{equation}\label{eqq}\dot{g}_i=T_{\overline{e}_i}L_{g_i}(u_i),\quad g_i(0)=g_0^{i},\end{equation} and the set $\mathcal{E}\subset\mathcal{V}\times\mathcal{V}$ denotes the set of edges of the graph, $s:=\frac{r(r-1)}{2}$ symmetric binary relations that link two agents, where $g_i(\cdot)\in C^{1}([0,T],G_i)$, $T\in\mathbb{R}$ fixed and $u_i$ is a curve on the Lie algebra $\mathfrak{g}_i$ of $G_i$. Alternatively, the left-invariant control system \eqref{eqq} can be written as $u_i(t)=T_{g_i}L_{g_i^{-1}}\dot{g}_i$, where for each $i$, the $m$-tuple of control inputs $u_i=[u^{1}_i\ldots u^{m}_i]^{T}$ take values in $\mathbb{R}^{m}$. 


The set $\mathcal{C}$ is given by $s$ (holonomic) formation constraints indexed by the edges set $c_{\mathcal{E}}=\{\phi_{ij}\}_{e_{ij}}$ with $e_{ij}=(v_i,v_j)\in\mathcal{E}$. For each edge $e_{ij}$, $\phi_{ij}$ is a function on $G_i\times G_j$ defining the formation constraint between agents $i$ and $j$. The constraint is enforced if and only if $\phi_{ij}(g_i,g_j)=0$.

If for each $i$, $\mathfrak{g}_i= \hbox{span}\{e_{1}^{i},\ldots,e_{m}^{i},e_{m+1}^{i},\ldots,e_{n}^{i}\}$, with $i=1,\ldots,r$, then $u_i$ is given by $\displaystyle{
u_i(t) = \sum_{w=1}^{m}u^{w}_i(t)e_{w}^{i}}$. Therefore  \eqref{eqq} is given by the \textit{drift-free kinematic left invariant control system} \begin{equation}\label{lics.1}\dot{g}_i=g_i\sum_{w=1}^{m}u^{w}_i(t)e_{w}^i.\end{equation}

Left-invariant control systems \eqref{eqq} provide a general framework for a class of systems that includes control design for spacecrafts and UAV motion. In general, the configuration space for these systems is globally described by a matrix Lie group making \eqref{eqq} a natural model for the controlled system. The Lie group framework gives rise to coordinate-free expressions for the dynamics describing the behavior of the system. When systems on Lie groups are left invariant, there is a globalization of  solutions. That is, even if we exploit local charts to make small maneuvers, working in a framework of Lie groups allow us to move all over the configuration space without reformulating the controls. This is because the absolute position of the system can always be described as if it were the identity in the Lie group.



\subsection{Problem formulation}

Next, we are going to define an optimal control problem for the left-invariant multi-agent control system \eqref{lics.1}. 

Along this work, we will denote $\displaystyle{G=\Pi_{i=1}^{r}G_{i}}$ and $\mathfrak{g}=\Pi_{i=1}^{r}\mathfrak{g}_{i}$, where the Lie algebra structure of $\mathfrak{g}$ is given by $[\xi_1,\xi_2]=([\xi_1^1,\xi_2^1],\ldots,[\xi_1^r,\xi_2^r])\in\mathfrak{g}$ with $\xi_1=(\xi_1^1,\ldots,\xi_1^r)\in\mathfrak{g}$ and $\xi_2=(\xi_2^1,\ldots,\xi_2^r)\in\mathfrak{g}.$
We also denote $\pi_{i}:G\to G_i$, $\tau_{i}:\mathfrak{g}\to\mathfrak{g}_{i}$, $\chi_i:\mathfrak{g}^{*}\to\mathfrak{g}_{i}^{*}$, $\alpha_i:T^{*}G\to T^{*}G_i$ and $\beta_i:TG\to TG_i$ the canonical projections from $G$, $\mathfrak{g}$, $\mathfrak{g}^{*}$, $T^{*}G$ and $TG$, respectively, over its $i^{th}$-factor. 

 We assume that each agent $i$ occupies a disk of radius $\overline{r}$ on $G$. The quantity $\overline{r}$ is chosen to be small enough so that it is possible to pack $r$ disks of radius $\overline{r}$ on $G$. 



We want to find necessary conditions for the existence of normal extremals in an energy-minimum problem for the left-invariant multi-agent control system \eqref{lics.1} where along their trajectory from a prescribed absolute initial state to a prescribed final absolute state not only minimize the energy of the complete networked system to achieve the desired final position, but also ensure that agents avoids collisions with each other in the workspace. This task can be done by introducing $s:=\frac{r(r-1)}{2}$ potential functions corresponding to fictitious forces into the cost functional for the optimal control problem, which are induced by the formation constraints.

Collision avoidance between the agents $i$ and $j$ is achieved by introducing the potential function $V_{ij}:G_i\times G_j\to\mathbb{R}$ with $i=1,\ldots,n$; $j=1,\ldots,s$ and $i\neq j$. We assume that $V_{ij}$ may not be a $(G_i\times G_j)$-invariant function for each $e_{ij}\in\mathcal{E}$ (i.e., $V_{ij}(g_ih_i,g_jh_j)\neq V_{ij}(g_i,g_j)$ for some $(g_i,g_j)$ and $(h_i,h_j)$ in $G_i\times G_j$) and it is sufficiently regular for all $e_{ij}\in\mathcal{E}$.


\textbf{Problem:} Consider the following optimal control problem: find the absolute configurations $g(t)=(g_1(t),\ldots,g_r(t))$ and control inputs $u(t)=(u_1(t),\ldots,u_r(t))$ minimizing the cost functional \begin{equation}\label{ocp}\min_{(g(\cdot), u(\cdot))}\sum_{i=1}^{r}\int_0^T[C_i(g_i(t),u_i(t))+\sum_{j>i}V_{ij}(g_i(t),g_j(t))]dt\end{equation} subject to $\dot{g}_i(t) = T_{\bar{e}_i}L_{g_i(t)}(u_i(t))$ and boundary values 
$g(0)=(g_1(0),\ldots,g_r(0)) =:(g^{0}_{1},\ldots,g^{0}_{r})$, $g(T)=(g_1(T),\ldots,g_r(T))=:(g^{T}_{1},\ldots,g^{T}_r)$, where $u(t)=(u_1(t),\ldots,u_r(t))\in\mathfrak{g}$ and each cost function $C_{i}:G_i\times\mathfrak{g}_{i}\to\mathbb{R}$ is smooth and $G_i$-invariant for each $i$ under the left action of $G_i$ on $G_i\times\mathfrak{g}_i$ given by $\rho_{g_i}: G_i\times\mathfrak{g}_i\to G_i\times\mathfrak{g}_i$, $\rho_{g_i}(h_i,u_i)=(L_{g_i}h_i,u_i)$. The cost functions $C_i$ are not related to collision avoidance between agents but only to the energy minimization performance.

\begin{remark}
For a given element $g_i\in G_{i}$, we denote $\rho_{g_i}:G_i\times\mathfrak{g}_i\to G_i\times\mathfrak{g}_i$, the action $\rho_{g_i}(h_i,u_i)=(L_{g_i}h_i,u_i)$, and when there is no confusion we denote $\rho:G\times (G\times\mathfrak{g})\to G\times\mathfrak{g}$ the vector valued left action of $G$ on $G\times\mathfrak{g}$ given by $\rho(g,h,u)=(L_gh,u)$ for $g,h\in G$, $u\in\mathfrak{g}$.\hfill$\diamond$
\end{remark}

\begin{remark}
The restriction in the index $j>i$ in the second sum in \eqref{ocp} is to do not count twice the quantity of functions $V_{ij}$ to avoid collision among agents (note that $V_{ij}=V_{ji}$) and make more simple the exposition. Other possibility might be to write in \eqref{ocp} the factor $\displaystyle{\frac{1}{2}\sum_{j=1,j\neq i}^{s}V_{ij}(g_i,g_j)}$ for $i=1,\ldots,r$. \hfill$\diamond$

\end{remark}



\section{Necessary conditions for the existence of normal extremals}

As in \cite{BlCoGuOh}, and \cite{KoonMarsden} for the single agent case, the optimal control problem can be solved as a constrained variational problem by introducing the Lagrange multipliers $\lambda_{g_i}=T^{*}_{g_i}L_{g_i^{-1}}(\lambda_i(t))\in T^{*}_{g_i}G_i$ with $\lambda_i\in C^{1}([0,T],\mathfrak{g}^{*})$ into the cost functional. Let $\mathfrak{g}^{*}=\hbox{span}\{e_i^{1},\ldots,e^{m}_i,e^{m+1}_i,\ldots,e^{n}\}$, then $\displaystyle{\lambda_{i}(t)=\sum_{w=m+1}^{n}\lambda_{i}(t)e_{i}^{w}}$, where the $(n-m)$-tuple of Lagrange multipliers  for the agent $i$, $\lambda_{i}=[\lambda^{m+1}_i,\ldots,\lambda^n_i]^{T}$ takes values in $\mathbb{R}^{(n-m)}$ for each $i$.

We define the set of admissible trajectories $\mathcal{A}\subset G\times\mathfrak{g}\times T^{*}G$ by \begin{align*}\mathcal{A}=\{(g,u,\lambda_g)\in G\times\mathfrak{g}&\times T^{*}G|\forall i=1,\ldots,r\\&\hbox{ and }j>i, \pi_i(g)-\pi_j(g)>\bar{r}\}\end{align*}
 and consider the extended Lagrangian $\mathcal{L}:\mathcal{A}\to\mathbb{R}$ given by \makeatletter 
 \def\@eqnnum{{\normalsize \normalcolor (\theequation)}} 
  \makeatother { \small
\begin{align*}\mathcal{L}(g,u,\lambda_g)=&\sum_{i=1}^{r}(C_{i}(\pi_{i}(g(t)),\tau_{i}(u(t)))\\&+\langle\alpha_{i}(\lambda_g(t)),\beta_{i}(T_{\overline{e}}L_{g}u(t))\rangle\\&+\sum_{j>i}V_{ij}(\pi_{i}(g(t)),\pi_j(g(t)))).\end{align*}}

The following result gives rise to necessary conditions for the existence of reduced normal extremals in the optimal control problem.
\begin{theorem}\label{theorem}
If the cost functions $C_{i}:G_i\times\mathfrak{g}_i\to\mathbb{R}$ are $G_i$-invariant under the left action of $G_i$ on $G_i\times\mathfrak{g}_i$ given by $\rho_{g_i}: G_i\times\mathfrak{g}_i\to G_i\times\mathfrak{g}_i$, $\rho_{g_i}(h_i,u_i)=(L_{g_i}h_i,u_i)$ for all $i=1,\ldots,r$ among the set $\mathcal{A}$ of admissible trajectories, the normal extremals for the optimal control problem \eqref{ocp} satisfy the Euler-Poincar\'e type equations 

\begin{align}
0=&\frac{d}{dt}\left(\frac{\partial C_i}{\partial u_i}+\lambda_i\right)-\hbox{ad}^{*}_{u_i}\left(\frac{\partial C_i}{\partial u_i}+\lambda_i\right)\label{eqep1}\\
&-\sum_{j>i}\left(T_{\overline{e}_{i}}^{*}L_{g_i}\left(\frac{\partial V_{ij}}{\partial g_i}\right)+T_{\overline{e}_j}^{*}L_{g_j}\left(\sum_{k=1}^{j-1}\frac{\partial V_{kj}}{\partial g_j}\right)\right)\nonumber
\end{align} together with $\dot{g}_{i}=g_iu_i$.
\end{theorem}
\begin{remark}\label{remark}
Note that equations  \eqref{eqep1} looks like they are not Euler-Poincar\'e equations due to the $g$-dependency, but we refer to them as Euler-Poincar\'e \textit{type} equations, with some abuse, (also known in the literature as trivialized Euler-Lagrange equations \cite{CoMdD}, \cite{Holmbook}) because they comes from a procedure of reduction by symmetries and if one considers the order in the way these equations must be solved, the first equation is an equation on the dual of the Lie algebra of $\mathfrak{g}$. In practice, these equations should be applied in backward sense with respect to the index $i$, as illustrated below.

For instance, in the case of three agents (the one for which we apply the reduced equations in Section \ref{section5}), $s=3$. We start by computing the equations when $i=3$. In this case $j=0$ and therefore we only need to solve the equation corresponding to the first two factors in \eqref{eqep1} since the other part disappears for $j=0$. Solving such an equation (we shown how such equation must be computed in Proposition \ref{proposition}) we obtain $u_3$ and $\lambda_3$, and therefore using $u_3$ and the reconstruction equation $\dot{g}_3=g_3u_3$ we obtain $g_3$. Therefore, after this first step we have the evolution of  $g_3, u_3$ and $\lambda_3$.

Next, we compute the equations when $i=2$. In this case $j=3$. Solving \eqref{eqep1} coupled with  
$\dot{g}_{2}=g_2u_2$, using the the configuration found for the last agent $g_3$ in the previous step, one is not able to solve the resulting equations and obtain the evolution for $\lambda_2$, $u_2$ and $g_2$ due to the interconnection term among the agents in $V_{13}$, although we have the information to deal with $V_{23}$ from the previous stage. Then, we must couple the equations with the case $i=1$, $j=2,3$. That is, at this stage we need to go to step $i=1$.

Solving such a system coupled also with $\dot{g}_{1}=g_1u_1$ we are able to find the evolution of $\lambda_1$, $\lambda_2$, $u_1$, $u_2$, $g_1$ and $g_2$, giving rise to the motion of the agents in the workspace and the control inputs, satisfying the formation constraints and minimizing the cost functional \eqref{ocp}.\hfill$\diamond$
\end{remark}

Note that equations \eqref{eqep1} can not describe completely the time evolution of the controls and the Lagrange multipliers. Since they are two independent variables, we must have two equations in order to have a system of differential equations with a well defined solution. To solve such a difficulty we propose the following splitting of the equations: 
 
\begin{proposition}\label{proposition}
If each $\mathfrak{g}_i$ and $\mathfrak{g}_{j}$ admits a decomposition of the form $\mathfrak{g}_i = \mathfrak{r}_i\oplus\mathfrak{s}_i$ and $\mathfrak{g}_j = \mathfrak{r}_j\oplus\mathfrak{s}_j$, respectively, where $\mathfrak{r}_i = \hbox{span}\{e_{1}^i,\ldots,e_{m}^i\}$, $\mathfrak{s}_i = \hbox{span}\{e_{m+1}^i,\ldots,e_{n}^i\}$ in agreement with \eqref{eqq}
and such that $$[\mathfrak{s}_i,\mathfrak{s}_i]\subseteq\mathfrak{s}_i,\quad [\mathfrak{s}_i,\mathfrak{r}_i]\subseteq\mathfrak{r}_i,\quad [\mathfrak{r}_i,\mathfrak{r}_i]\subseteq\mathfrak{s}_i,$$ $$[\mathfrak{s}_j,\mathfrak{s}_j]\subseteq\mathfrak{s}_j,\quad [\mathfrak{s}_j,\mathfrak{r}_j]\subseteq\mathfrak{r}_j,\quad [\mathfrak{r}_j,\mathfrak{r}_j]\subseteq\mathfrak{s}_j,$$ for each $i=1,\dots,r$, with $j=1,\ldots,s$ where $s=\frac{r(r-1)}{2}$, then the time evolution of \eqref{eqep1} for each $i$ can be rewritten as
\begin{align*}
\frac{d}{dt}\frac{\partial C_i}{\partial u_i} =& \ad_{u_i}^{*}\lambda_i\bigg|_{\mathfrak{r}_i}+\sum_{j>i}T_{\overline{e}_j}^{*}L_{g_j}\left(\sum_{k=1}^{j-1}\frac{\partial V_{kj}}{\partial g_j}\right)\bigg|_{\mathfrak{r}_j}\nonumber\\
&+\sum_{j>i}T_{\overline{e}_{i}}^{*}L_{g_i}\left(\frac{\partial V_{ij}}{\partial g_i}\right)\bigg|_{\mathfrak{r}_i}\\
\dot{\lambda}_i =& \ad_{u_i}^{*}\frac{\partial C_i}{\partial u_i}\bigg|_{\mathfrak{s}_i}+\sum_{j>i}T_{\overline{e}_j}^{*}L_{g_j}\left(\sum_{k=1}^{j-1}\frac{\partial V_{kj}}{\partial g_j}\right)\bigg|_{\mathfrak{s}_j}\nonumber\\
&+\sum_{j>i}T_{\overline{e}_{i}}^{*}L_{g_i}\left(\frac{\partial V_{ij}}{\partial g_i}\right)\bigg|_{\mathfrak{s}_i}\nonumber\\
\end{align*}
\end{proposition}

\section{Optimal control of multiples unicycles}\label{section5}

As an application we consider an energy-minimum control problem for three unicycle type robots. A unicycle is a homogeneous disk rolling on a horizontal plane maintaining its vertical position (see, e.g. \cite{Bl}). The configuration of each unicycle at any given time is determined by the element $g_i\in\mathrm{SE}(2)\cong\mathbb{R}^{2}\times\mathrm{SO}(2)$, $i=1,2,3$ given by $\displaystyle{g_i = \begin{bmatrix}
    \cos\theta_i & -\sin\theta_i & \phantom{-}x_i\\
    \sin\theta_i & \phantom{-}\cos\theta_i & \phantom{-}y_i\\
    0 & \phantom{-}0 & \phantom{-}1
  \end{bmatrix},\nonumber}$ where $(x_i, y_i)\in\mathbb{R}^{2}$ represents the point of contact of each wheel with the ground and $\theta_i\in SO(2)$ represents the angular orientation of each unicycle. We denote $u_i=(u_i^1,u_i^2)$ with $i=1,2,3$. The kinematic equations for the multi-agent system are
\begin{equation}\label{unicycle}
\dot{x}_i = u^{2}_i\cos\theta_i,\quad
\dot{y}_i = u^{2}_i\sin\theta_i,\quad
\dot{\theta}_i= u^{1}_i\quad i=1,2,3.
\end{equation}

%

Equations \eqref{unicycle} gives rise to a left-invariant control system on $\mathrm{SE}(2)\times SE(2)\times SE(2)$ (note that the first two equations are equivalent to the nonholonomic constraints $\dot{x}_i\sin\theta_i-\dot{y}_i\cos\theta_i= 0$). As a left invariant control system on $SE(2)\times SE(2)\times SE(2)$ the unicycle equations take the form $$\dot{g}_i=g_i(e_1^{i}u_i^{1}+e_2^{i}u_i^{2})$$ describing all directions of allowable motion, where the elements of the basis of $\mathfrak{g}_i=\mathfrak{se}(2)_i$ are 
\makeatletter 
 \def\@eqnnum{{\normalsize \normalcolor (\theequation)}} 
  \makeatother

{ \small
\begin{equation*}
e_{1}^i = \begin{bmatrix}
        0 & -1 & \phantom{-}0\\
        1 & \phantom{-}0 & \phantom{-}0\\
        0 & \phantom{-}0 & \phantom{-}0
        \end{bmatrix}
e_{2}^i = \begin{bmatrix}
        0 & \phantom{-}0 & \phantom{-}1\\
        0 & \phantom{-}0 & \phantom{-}0\\
        0 & \phantom{-}0 & \phantom{-}0
        \end{bmatrix} e^i_{3} = \begin{bmatrix}
        0 & \phantom{-}0 & \phantom{-}0\\
        0 & \phantom{-}0 & \phantom{-}1\\
        0 & \phantom{-}0 & \phantom{-}0
        \end{bmatrix}\end{equation*}}which satisfy
$[e_{1}^i,e_{2}^i] = e_{3}^i, \hspace{5pt} [e_{2}^i,e_{3}^i] = 0_{3\times 3}, \hspace{5pt} [e_{3}^i,e_{1}^i] = e_{2}^i$. 
Using the dual pairing, where $\langle\alpha_i,\xi_i\rangle:=\hbox{tr}(\alpha_i\xi_i)$, for any $\xi_i\in\mathfrak{se}(2)_i$ and $\alpha_i\in\mathfrak{se}(2)_i^{*}$, the elements of the basis of $\mathfrak{se}(2)_i^{*}$ are given by
\makeatletter 
 \def\@eqnnum{{\normalsize \normalcolor (\theequation)}} 
  \makeatother

{ \footnotesize\begin{equation*}
e^{1}_i = \begin{bmatrix}
        \phantom{-}0 & \phantom{-}\dfrac{1}{2} & \phantom{-}0\\
        -\dfrac{1}{2} & \phantom{-}0 & \phantom{-}0\\[8pt]
        \phantom{-}0 & \phantom{-}0 & \phantom{-}0
        \end{bmatrix} e^{2}_i = \begin{bmatrix}
        0 & \phantom{-}0 & \phantom{-}0\\
        0 & \phantom{-}0 & \phantom{-}0\\
        1 & \phantom{-}0 & \phantom{-}0
        \end{bmatrix} e^{3}_i = \begin{bmatrix}
        0 & \phantom{-}0 & \phantom{-}0\\
        0 & \phantom{-}0 & \phantom{-}0\\
        0 & \phantom{-}1 & \phantom{-}0
        \end{bmatrix}.\end{equation*}}
Here, $\mathfrak{r}_{i}=\{e_1^i,e_2^i\}$, $\mathfrak{s}_{i}=\{e_3^i\}$, $\mathfrak{se}(2)_i=\mathfrak{r}_{i}\oplus\mathfrak{s}_{i}$ and fulfill the hypothesis of Proposition \ref{proposition}. Also note that $\mathfrak{se}(2)$ is not a semi-simple Lie algebra but it satisfies the Lie bracket relations to decompose the dynamics as in Proposition \ref{proposition}.

The formation is completely specified by the (holonomic) constraints $\phi_{ij}:SE(2)_{i}\times SE(2)_{j}\to\mathbb{R}$, $j>i$, $i,j=1,2,3$ (i.e., $\phi_{12},\phi_{13},\phi_{23}$, $s=3$) determined by a prescribed distance $d_{ij}$ among the center of masses of all agent at any time. The constraint for the edge $e_{ij}$ is given by  $\phi_{ij}(g_i,g_j)=||\psi(g_j)g_i||^2_{F}-\tilde{d}_{ij}$ where $||\cdot||_{F}$ is the Frobenius norm, $||A||_{F}=\hbox{tr}\left(A^{T}A\right)^{1/2}$, $\tilde{d}_{ij}=d_{ij}^2+3$ and $\psi:SE(2)\to SE(2)$ is the smooth map defined as $\psi(g)=\bar{g}$ where $\bar{g}=\left[ {\begin{array}{ccc}
   1 & 0 & -x \\
   0 & 1 & -y \\
   0 & 0 &1\\
  \end{array} } \right]\in SE(2).$ 
  
  It is straightforward to corrobarate that the constraint $\phi_{ij}(g_i,g_j)=0$ on absolute configurations on the Lie group $SE(2)_{i}\times SE(2)_{j}$, is equivalent to the constrain in the relative configurations $(x_i-x_j)^2+(y_i-y_j)^2-d_{ij}^2$.


 
 
Consider the functions $V_{ij}:\mathbb{R}^{2}\times\mathbb{R}^{2}\to\mathbb{R}$ given by $\displaystyle{V_{ij}((x_i,y_i),(x_j,y_j)) = \frac{\sigma_j}{2((x_i-x_j)^{2}+(y_i-y_j)^{2}-d_{ij}^{2})}}$, where $\sigma_j\in\mathbb{R}^{+}$. The inner product on $\mathfrak{se}(2)_{i}$ is given by $\langle\langle\xi_i,\xi_i\rangle\rangle=\hbox{tr}(\xi_i^{T}\xi_i)$, for any $\xi_i\in\mathfrak{se}(2)_i$ and hence, the norm $\|\xi_i\|_{\mathfrak{se}(2)}$ is given by $\|\xi_i\|_{\mathfrak{se}(2)}=\langle\langle\xi_i,\xi_i\rangle\rangle^{1/2} = \sqrt{\hbox{tr}(\xi_i^{T}\xi_i)}$, for any $\xi_i\in\mathfrak{se}(2)_i$. 

The potential functions defined on the Lie group are given by $V_{ij}: \mathrm{SE}(2)_i\times\mathrm{SE}(2)_j\to\mathbb{R}$,  
$$V_{ij}(g_i,g_j) = \frac{\sigma_j}{2(\|\psi(g_j)g_i\|_{F}^{2}-\tilde{d}_{ij})}$$ and the cost functions are given by $C_i(g_i,u_i)= \frac{1}{2}\langle u_i,u_i\rangle$
where $u_i= u^{1}_ie_{1}^i+u^{2}_ie_{2}^i$. 

The potential functions $V_{ij}$ are $SO(2)\times SO(2)$-invariants but not $SE(2)\times SE(2)$-invariants as the cost functions, and therefore the $SE(2)$-symmetry (or invariance) in the dynamics of the agents is broken in the optimal control problem due to the incorporation of the formation constraints into the cost functional.

The Lagrangian for the optimal control problem $\mathcal{L}:(SE(2))^{3}\times(\mathfrak{se}(2))^{3}\times(T^{*}SE(2))^{3}\to\mathbb{R}$ is given by \begin{align*}\mathcal{L}(g_i,u_i,\lambda_{g_i})=&\sum_{i=1}^{3}\frac{1}{2}\langle u_i,u_i\rangle+\langle T^{*}_{g_i}L_{g_i^{-1}}(\lambda_i),g_iu_i\rangle\\&+\sum_{j>i} \frac{\sigma_j}{2(\|\psi(g_j)g_i\|_{F}^{2}-\tilde{d}_{ij})}\end{align*} where $\lambda_i=\lambda_i^3e_i^3$. The reduced Lagrangian $\ell:(SE(2))^{3}\times(\mathfrak{se}(2))^{3}\times(\mathfrak{se}(2)^{*})^{3}\to\mathbb{R}$ is given by \begin{align*}\ell(g_i,u_i,\lambda_i)=&\sum_{i=1}^{3}\frac{1}{2}\langle u_i,u_i\rangle+\langle \lambda_i,u_i\rangle\\&+\sum_{j>i} \frac{\sigma_j}{2(\|\psi(g_j)g_i\|_{F}^{2}-\tilde{d}_{ij})}\end{align*}

Note that $\frac{\partial C_i}{\partial u_i}=u_i$, \begin{align}
\ad_{u_i}^{*}\lambda_i =& \begin{bmatrix} 
                      0 & -\dfrac{u^{2}_i\lambda^{3}_i}{2} & \phantom{-}0\\
										  \dfrac{u^{2}_i\lambda^{3}_i}{2} & \phantom{-}0 & \phantom{-}0\\[8pt] 
										  u^{1}_i\lambda^{3}_i & \phantom{-}0 & \phantom{-}0
										  \end{bmatrix},\nonumber\\ \ad_{u_i}^{*}u_i =& \begin{bmatrix} 
                0 & \phantom{-}0 & \phantom{-}0\\
							  0 & \phantom{-}0 & \phantom{-}0\\
							  0 & -u^{1}_iu^{2}_i & \phantom{-}0
							  \end{bmatrix} \nonumber
\end{align}and
 $$T_{\overline{e}_i}L_{g_i}\left(\frac{\partial \overline{V}_{ij}}{\partial g_i}\right)\Big{|}_{\mathfrak{r}_i} = 
-\frac{\sigma_j([g_i^{T}\psi(g_j)^{T}\psi(g_j)g_i]_{13})}{2(\|\psi(g_j)g_i\|^{2}_{F}-\tilde{d}_{ij})^{2}}e_2^{i}$$								
$$T_{\overline{e}_i}L_{g_i}\left(\frac{\partial \overline{V}_{ij}}{\partial g_i}\right)\Big{|}_{\mathfrak{s}_i }= 
-\frac{\sigma_j([g_i^{T}\psi(g_j)^{T}\psi(g_j)g_i]_{23})}{2(\|\psi(g_j)g_i\|^{2}_{F}-\tilde{d}_{ij})^{2}}e_3^i$$ where the subindex $13$ and $23$ stands for the entry $13$ and $23$ of the matrix $g_{ij}:=g_i^{T}\psi(g_j)^{T}\psi(g_j)g_i$ and the partial derivatives with respect to $g_j$ are computed in a similar fashion where the unique difference is a sign factor (instead of $-\sigma_j$ must be $+\sigma_j$).  

By Remark \ref{remark}  and using Proposition \ref{proposition} the Euler-Poincar\'e equations giving rise to the reduced necessary conditions for the existence of extremals in the optimal control problem for $i,j=1,2,3$, are computed in backward order in the index $i$ as follows.

When $i=3$, then $j=0$ and therefore the equations obtained by using Proposition \ref{proposition} are 
 
$$\dot{u}^{1}_3 =-\frac{u_3^2\lambda^3_3}{2},\quad\dot{u}^2_3=u^1_3\lambda_3^3,\quad\dot{\lambda}_{3}^{3}=-u_3^1u_3^2$$
where the subindex $3$ in all the equations stands for the agent number $3$.

The second step consist on solving the system of equations for $i=2$. In this case $j=3$ and the resulting differential equations are 
\begin{align*}
\dot{g_2}=&g_2(e_1^2u_2^1+e_2^2u_2^2)\\
\dot{u}_{2}^1=&-\frac{u_2^2\lambda_2^3}{2}\\
\dot{u}_{2}^{2}=&u_2^1\lambda_2^3-\frac{\sigma_3(g_{23})_{13}}{(\|\psi(g_3)g_2\|^{2}_{F}-\tilde{d}_{23})^{2}}+\frac{\sigma_3(g_{13})_{13}}{(\|\psi(g_3)g_1\|^{2}_{F}-\tilde{d}_{13})^{2}}\\
&+\frac{\sigma_3(g_{23})_{13}}{(\|\psi(g_3)g_2\|^{2}_{F}-\tilde{d}_{23})^{2}}\\
\dot{\lambda}_{2}^{3}=&-u_2^1u_2^2-\frac{\sigma_3(g_{23})_{23}}{(\|\psi(g_3)g_2\|^{2}_{F}-\tilde{d}_{23})^{2}}+\frac{\sigma_3(g_{13})_{23}}{(\|\psi(g_3)g_1\|^{2}_{F}-\tilde{d}_{13})^{2}}\\
&+\frac{\sigma_3(g_{23})_{23}}{(\|\psi(g_3)g_2\|^{2}_{F}-\tilde{d}_{23})^{2}}.
\end{align*} Such a system of equations can not provide the evolution of the configurations, controls and Lagrange multipliers $g_2$, $u_2$ and $\lambda_2$ respectively by the interconnection term between $g_1$ and $g_3$, and we must to couple the previous equations with the ones for $i=1$.
 
Finally, computing the equations provided by Proposition \ref{proposition} when $i=1$, $j=2,3$ we have
 
 \begin{align*}
\dot{g_1}&=g_1(e_1^1u_1^1+e_2^1u_1^2)\\
\dot{u}_{1}^1&=-\frac{u_1^2\lambda_1^3}{2}\\
\dot{u}_{1}^{2}&=u_1^1\lambda_1^3-\frac{\sigma_3(g_{23})_{13}}{(\|\psi(g_3)g_2\|^{2}_{F}-\tilde{d}_{23})^{2}}\\
\dot{\lambda}_{1}^{3}&=-u_1^1u_1^2 +\frac{\sigma_3(g_{23})_{23}}{(\|\psi(g_3)g_2\|^{2}_{F}-\tilde{d}_{23})^{2}}.
\end{align*}
Solving the previous system of equations for $i=1$ and $i=2$ together, we can obtain the solutions for all absolute configurations, controls and Lagrange multipliers that permits to move agents from a prescribed initial states to a prescribed final state minimizing the total energy in the multi-agent formation and avoiding collision between the three unicycles.

\section{Conclusions and Future Research}
We studied reduced necessary conditions for the existence of normal extremals in an optimal control problem for a multi-agent formation, defined by kinematic drift-free left invariant control systems with formation constraints breaking the symmetry of the formation graph. Such conditions were obtained through a non-trivial extension of variational principles for a single agent whose dynamics evolves on a Lie group. Such a variational principle does not allow to describe the conditions of extremals for the optimal control problem as a solvable and well defined system of first order differential equations for the evolution of the agents while at the same time minimizing the energy and satisfying the formation constraints. We proposed a splitting into the dynamics to make the system solvable.

As an application we studied the optimal control of multiples unicycles and we characterized some integrability properties to consider when one solves the corresponding system of equations to obtain the optimal trajectories satisfying the formation constraints.

In a future work we will explore the incorporation of decentralized navigation functions \cite{dimos}, \cite{TK}, \cite{chris} into the cost functional to avoid simultaneously collision between agents and static obstacles in the workspace as well as agents with limiting sensing \cite{Freeman}. 



\begin{thebibliography}{20}
\providecommand{\newblock}{\relax}
\bibitem{balch} T. Balch and R. Arkin. Behavior-based formation control for multirobot systems. IEEE Trans. Robot. Autom. Vol 14 (6), 926-939, 1998.
\bibitem{Bl}
A.~M. Bloch. Nonholonomic mechanics and control. Series IAM. New York: Springer-Verlag, vol.~24. 2nd Edition 2015.
\bibitem{BlCoGuOh} A. Bloch, L. Colombo, R. Gupta and T. Ohsawa. Optimal control problems with symmetry breaking cost functions. SIAM J. Applied Algebra and Geometry, 1(1), 626-646.  
\bibitem{CoMdD} L. Colombo and D. Mart\'in de Diego. Higher-order variational problems
on Lie groups and optimal control applications, J. Geom. Mech., vol. 6 (4), 451-478, 2014.
\bibitem{dimos} D. Dimarogonas, S. Loizou, K. Kyriakopoulos and M. Zavlanos. A feedback stabilization and collision avoidance scheme for multiple independent non-point agents. Automatica, Vol. 42 (2), 229-243, 2006.
\bibitem{Freeman} R. Freeman, P. Yang and K. Lynch.  Distributed estimation and control of swarm formation statistics.  American Control Conf., 749-755, 2006.

\bibitem{Holmbook} D. Holm, T. Schmah and C Stoica. Geometric mechanics and symmetry:  from finite to infinite dimensions. Oxford University Press, 2009
\bibitem{JK}E. Justh and P. Krishnaprasad. Equilibria and steering laws for planar formations. Systems $\&$ Control Letters 52, 25-38, 2004.
\bibitem{KoonMarsden} W. Koon and J. Marsden. Optimal control for holonomic and nonholonomic mechanical systems with symmetry and Lagrangian reduction. SIAM Journal on Control and Optimization, 35 (3), 901-929, 1997.

\bibitem{Krishnaprasad} P. Krishnaprasad. Optimal control and Poisson reduction. Technical Report T.R. 93-87, University of Maryland, 1993.

\bibitem{Leonard1} N. Leonard and P. Krishnaprasad, Motion control of drift-free,
left-invariant systems on lie groups, IEEE Transactions on Automatic Control,
vol. 40 (9), 1539-1554, 1995.
\bibitem{Leonard2} N. Leonard and  E. Fiorelli. Virtual leaders, artificial potentials and coordinated control of groups. IEEE Conference on Decision and Control, 2968-2973, 2001.



\bibitem{MMbook} M. Mesbahi and M. Egerstedt. Graph theoretic methods in multiagent networks. Princeton University Press, 2010.
\bibitem{Sastrybook} R. Murray, Z. Li and S. Sastry. A mathematical introduction to robotic manipulation. CRC Press, 1994.
\bibitem{pedro1} P. O. Pereira and D. V. Dimarogonas. Family of controllers for attitude synchronization on the sphere. Automatica, vol. 75, pp. 271-281, 2017.
\bibitem{pedro2} P. O. Pereira, R. Cunha, D. Cabecinhas, C. Silvestre, P. Oliveira. Leader following trajectory planning: A trailer-like approach. Automatica, vol. 75, pp. 77-87, 2017.
\bibitem{tabuada} P. Tabuada, J. Pappas and P. Lima. Motion feasibility dor multi-agent formation. IEEE Transactions on robotics. vol 21 (3), 387-392, 2005.
\bibitem{TK} H.  Tanner  and  A.  Kumar.  Formation  stabilization  of  multiple agents using decentralized navigation functions. Robotics: Science and Systems, MIT Press, 49- 56, 2005.
 \bibitem{chris} C. K. Verginis, A. Nikou and D. V. Dimarogonas. On the Position and Orientation Based Formation Control of Multiple Rigid Bodies with Collision Avoidance and Connectivity Maintenance. 56th IEEE Conference on Decision and Control (CDC), 2017, pp 411-416.



\end{thebibliography}
\end{document}